\documentclass{amsart}
\usepackage{amscd,amssymb,amsxtra}

\theoremstyle{plain}
\newtheorem{Thm}{Theorem}[section]
\newtheorem{Prop}[Thm]{Proposition}
\newtheorem{Lem}[Thm]{Lemma}
\newtheorem{Cor}[Thm]{Corollary}
\newtheorem{Claim}[Thm]{Claim}

\theoremstyle{definition}
\newtheorem{Def}[Thm]{Definition}
\newtheorem{Question}[Thm]{Question}

\theoremstyle{remark}

\newtheorem{Remark}[Thm]{Remark}
\newtheorem{Example}[Thm]{Example}
\newtheorem{Application}[Thm]{Application}

\numberwithin{equation}{section}

\newcommand{\supp}{\operatorname{supp}}
\newcommand{\res}{\upharpoonright}

\newcommand{\Aut}{\operatorname{Aut}}
\newcommand{\Inn}{\operatorname{Inn}}

\newcommand{\cf}{\operatorname{cf}}

\newcommand{\Sym}{\operatorname{Sym}}
\newcommand{\Alt}{\operatorname{Alt}}
\newcommand{\Wr}{\text{wr}}

\begin{document}

%Topmatter
\title[The automorphism tower problem revisited]
{The automorphism tower problem revisited}
%Author info
\author{Winfried Just}
\address
{Department of Mathematics \\
Ohio University \\
Athens \\
Ohio 45701}
\thanks{The research of the first author was done during a visit
to Rutgers University in September 1997, which was partially 
supported by NSF Grant DMS-9704477.}
\author{Saharon Shelah}
\address
{Mathematics Department \\
The Hebrew University \\
Jerusalem \\
Israel}
\thanks{The research of the second author was partially supported
by the U.S.-Israel Binational Science Foundation. This paper is
number 654 in the cumulative list of the second author's
publications.}
\author{Simon Thomas}
\address
{Mathematics Department \\
Rutgers University \\
110 Frelinghuysen Road \\
Piscataway \\
New Jersey 08854-8019}
\thanks{The research of the third author was partially supported
by NSF Grants.}

%End topmatter

\begin{abstract}
It is well-known that the automorphism towers of infinite 
centreless groups of cardinality $\kappa$ terminate in
less than $\left( 2^{\kappa} \right)^{+}$ steps. But an
easy counting argument shows that 
$\left( 2^{\kappa} \right)^{+}$ is not the best possible
bound. However, in this paper, we will show that it is
impossible to find an explicit better bound using $ZFC$.
\end{abstract}

\maketitle

\section{Introduction} \label{S:intro}
If $G$ is a centreless group,
then there is a natural embedding  of $G$ into its automorphism
group $\Aut G$, obtained by sending each $g\in G$ to the
corresponding inner automorphism $i_{g} \in\Aut G$. In this
paper, we will always work with the left action of $\Aut G$
on $G$.  Thus $i_{g}(x) = g x  g^{-1}$ for all $x \in G$.  
If $\pi \in \Aut G$ and $g \in
G$, then an easy calculation shows that
$\pi i_{g} \pi^{-1} = i_{\pi(g)}$.  Hence the group of inner
automorphisms $\Inn G$ is a normal subgroup of $\Aut G$; and
$C_{\Aut G} (\Inn G) = 1$.  In particular, $\Aut G$ is also a centreless group.  
This enables us to define the automorphism tower of $G$ to be the
ascending chain of groups
\[
G = G_{0} \trianglelefteq G_{1} \trianglelefteq G_{2} \trianglelefteq  \dots
G_{\alpha} \trianglelefteq G_{\alpha +1} \trianglelefteq \dots
\]
such that for each ordinal $\alpha$
\begin{enumerate}
\item[(a)] $G_{\alpha +1} = \Aut G_{\alpha}$; and
\item[(b)] if $\alpha$ is a limit ordinal, then
$G_{\alpha} = \underset{\beta < \alpha}{\bigcup}G_{\beta}$.
\end{enumerate}
(At each successor step, we identify $G_{\alpha}$ with
$\Inn G_{\alpha}$ via the natural embedding.)

The automorphism tower is said to terminate if there exists an ordinal
$\alpha$ such that $G_{\alpha +1} = G_{\alpha}$.
This occurs if and only if there exists an ordinal $\alpha$ such that
$\Aut G_{\alpha}= \Inn G_{\alpha}$.  A
classical result of Wielandt \cite{w} says that if $G$ is finite, then the
automorphism tower terminates after finitely many steps. Wielandt's 
theorem fails for infinite centreless groups. For example, consider the
infinite dihedral group $D_{\infty}=\langle a,b \rangle$, where $a$ and $b$
are elements of order 2. Then
$D_{\infty} = \langle a \rangle \ast \langle b \rangle$ is the free product of 
its cyclic subgroups $\langle a \rangle$ and $\langle b \rangle$. It
follows that $D_{\infty}$ is a centreless group; and that
$D_{\infty}$ has an outer automorphism $\pi$ of order 2 which interchanges
the elements $a$ and $b$. It is easily shown that
$\Aut D_{\infty} = \langle \pi, i_{a} \rangle$. Thus $\Aut D_{\infty}$ is
also an infinite dihedral group, and so $\Aut D_{\infty} \simeq D_{\infty}$.
Hence for each $n \in \omega$, the $n^{th}$ group in the automorphism
tower of $D_{\infty}$ is isomorphic to $D_{\infty}$; and the automorphism tower 
of $D_{\infty}$ does not terminate after finitely many steps.

In the 1970s, a number of special cases of the automorphism tower problem
were solved. For example, Rae and Roseblade \cite{rr} proved that the automorphism
tower of a centreless \v{C}ernikov group terminates after finitely
many steps; and Hulse \cite{hu} proved that the automorphism tower
of a centreless polycyclic group terminates in countably many
steps. But the problem was not solved in full generality until 1984,
when Thomas \cite{t1} showed that the automorphism tower of an arbitrary
centreless group eventually terminates; and that for each ordinal $\alpha$,
there exists a group whose automorphism tower terminates in exactly
$\alpha$ steps.

\begin{Def} \label{D:height}
If $G$ is a centreless group, then the {\em height\/} $\tau(G)$
of the automorphism tower of $G$ 
is the least ordinal $\alpha$ such that $G_{\alpha +1} = G_{\alpha}$.
\end{Def}

This raises the question of finding bounds for $\tau(G)$ in terms
of the cardinality of $G$. In his original paper \cite{t1},
Thomas proved that if $G$ is an infinite centreless group of 
cardinality $\kappa$, then $\tau(G) \leq \left( 2^{\kappa} \right)^{+}$.
Soon afterwards, Thomas and Felgner independently noticed that
an easy application of Fodor's Lemma yielded the following
slightly better bound.

\begin{Thm}[Thomas \cite{t2}] \label{T:bound}
If $G$ is an infinite centreless group of cardinality $\kappa$, then
$\tau (G) < \left( 2^{\kappa} \right)^{+}$.
\end{Thm}
\begin{flushright}
$\square$
\end{flushright}

\begin{Def} \label{D:tor}
If $\kappa$ is an infinite cardinal, then $\tau_{\kappa}$
is the least ordinal such that $\tau(G) < \tau_{\kappa}$ for every
centreless group $G$ of cardinality $\kappa$.
\end{Def}

Since there are only $2^{\kappa}$ centreless groups of cardinality
$\kappa$ up to isomorphism, it follows that
$\tau_{\kappa} < \left( 2^{\kappa} \right)^{+}$. On the other hand,
Thomas \cite{t1} has shown that
for each ordinal $\alpha < \kappa^{+}$, there exists a centreless
group $G$ of cardinality $\kappa$ such that $\tau (G) = \alpha$.
Thus 
$\kappa^{+} \leq \tau_{\kappa} < \left( 2^{\kappa} \right)^{+}$.
It is natural to ask whether a better explicit bound on 
$\tau_{\kappa}$ can be proved in $ZFC$, preferably one which
does not involve cardinal exponentiation.

The proof of Theorem \ref{T:bound} is extremely simple, and uses only 
the most basic results in group theory, together with some
elementary properties of the infinite cardinal numbers. So it is
not surprising that Theorem \ref{T:bound} does not give the best possible
bound for $\tau_{\kappa}$. In contrast, the proof of Wielandt's
theorem is much deeper, and involves an intricate study of the
subnormal subgroups of a finite centreless group. The real question
behind the search for better explicit bounds for $\tau_{\kappa}$
is whether there exists a subtler, more informative, group-theoretic
proof of the automorphism tower theorem for infinite groups. The
main result of this paper says that no such bounds can be proved in
$ZFC$, and thus can be interpreted as saying that no such proof
exists. (It is perhaps worth mentioning that the proof of
Theorem \ref{T:bound} yields that the automorphism tower of a finite
centreless group terminates in {\em countably\/} many steps.
However, there does not seem to be an easy reduction from
countable to finite; and it appears that some form of Wielandt's
analysis is necessary.)

\begin{Thm} \label{T:best}
Let $V \vDash GCH$ and let $\kappa$, $\lambda  \in V$ be uncountable cardinals
such that $\kappa < \cf (\lambda)$. Let $\alpha$ be any ordinal such
that $\alpha < \lambda^{+}$. Then there exists a notion of forcing
$\mathbb{P}$, which preserves cofinalities and cardinalities, such
that the following statements are true in the corresponding generic
extension $V^{\mathbb{P}}$.
\begin{enumerate}
\item[(a)] $2^{\kappa} = \lambda$.
\item[(b)] There exists a centreless group $G$ of cardinality $\kappa$
such that $\tau (G) = \alpha$.
\end{enumerate}
\end{Thm}

Thus it is impossible to find better explicit bounds for $\tau_{\kappa}$
when $\kappa$ is an uncountable cardinal. However, our methods do
not enable us to deal with countable groups; and
it remains an open question whether or not
there exists a countable centreless group $G$ such that
$\tau (G) \geq \omega_{1}$.

Most of this paper will be concerned with the problem of constructing
centreless groups with extremely long automorphism towers.
Unfortunately it is usually very difficult to compute the
successive groups in an automorphism tower. We will get around this
difficulty by reducing it to the much easier computation of the
successive normalisers of a subgroup $H$ of a group $G$.

\begin{Def} \label{D:normaliser}
If $H$ is a subgroup of the group $G$, then the {\em normaliser tower\/}
of $H$ in $G$ is defined inductively as follows.
\begin{enumerate}
\item[(a)] $N_{0}(H) = H$.
\item[(b)] If $\alpha = \beta + 1$, then
$N_{\alpha}(H) = N_{G} \left( N_{\beta}(H) \right)$.
\item[(c)] If $\alpha$ is a limit ordinal, then
$N_{\alpha}(H) = \underset{\beta < \alpha}{\bigcup}N_{\beta}(H)$.
\end{enumerate}
It is sometimes necessary for the notation to include an explicit
reference to the ambient group $G$. In this case, we will write
$N_{\alpha}(H) = N_{\alpha}(H,G)$.
\end{Def}

As we will see in Section \ref{S:normaliser},
if $\alpha$ is any ordinal, then it is easy to construct examples of pairs
of groups, $H < G$, such that the normaliser tower of $H$ in $G$
terminates in exactly $\alpha$ steps. 
The following lemma, which was essentially proved in \cite{t1},
will enable us to convert  normaliser towers into 
corresponding automorphism towers. 

\begin{Lem} \label{L:psl}
Let $K$ be a field such that $|K| > 3$ and let $H$ be a subgroup of $\Aut K$.
Let
\[
G = PGL(2,K) \rtimes H \leqslant P \varGamma L(2,K) =
PGL(2,K) \rtimes \Aut K.
\]
Then $G$ is a centreless group; and for each $\alpha$,
$G_{\alpha} = PGL(2,K) \rtimes N_{\alpha}(H)$, where
$N_{\alpha}(H)$ is the $\alpha^{th}$ group in the normaliser
tower of $H$ in $\Aut K$. 
\end{Lem}
\begin{flushright}
$\square$
\end{flushright}

It is well-known that every group $G$ can be realised as the 
automorphism group of a suitable graph $\Gamma$. Thus the
following result implies that every group $G$ can also be
realised as the automorphism group of a suitable field $K$.

\begin{Lem}[Fried and Koll\'{a}r \cite{fk}] \label{L:field}
Let $\Gamma = \langle X, E \rangle$ be any graph. Then there exists a
field $K_{\Gamma}$ of cardinality $\max \{ \left| X \right|,  \omega \}$
which satisfies the following conditions.
\begin{enumerate}
\item[(a)] $X$ is an $\Aut 	K_{\Gamma}$-invariant subset of $K_{\Gamma}$.
\item[(b)] The restriction mapping, $\pi \mapsto \pi \res X$, is an 
isomorphism from $\Aut K_{\Gamma}$ onto $\Aut \Gamma$.
\end{enumerate}
\end{Lem}
\begin{flushright}
$\square$
\end{flushright}

Combining Lemmas \ref{L:psl} and \ref{L:field}, we obtain the following
reduction of our problem.

\begin{Lem} \label{L:reduction}
Suppose that there exists a graph $\Gamma$ of cardinality $\kappa$ and a
subgroup $H$ of $\Aut \Gamma$ such that
\begin{enumerate}
\item[(a)] $\left| H \right| \leq \kappa$; and
\item[(b)] the normaliser tower of $H$ in $\Aut \Gamma$ terminates in exactly 
$\alpha$ steps.
\end{enumerate}
Then there exists a centreless group $T$ of cardinality $\kappa$ such that
$\tau ( T) = \alpha$.
\end{Lem}

\begin{proof}
Let $K_{\Gamma}$ be the
corresponding field, which is given by Lemma \ref{L:field}, and let
$T = PGL(2,K_{\Gamma}) \rtimes H$. By Lemma \ref{L:psl},
$\tau (T) = \alpha$.
\end{proof}

\relax From now on, fix a regular uncountable cardinal $\kappa$ such that 
$\kappa^{< \kappa} = \kappa$ and an ordinal $\alpha$. Roughly
speaking, our strategy will be to
\begin{enumerate}
\item[(1)] first construct a pair of groups, $H < G$, such that 
$\left| H \right| \leq \kappa$ and the normaliser tower of $H$ in $G$ terminates
in $\alpha$ steps; and
\item[(2)] then attempt to find a cardinal-preserving notion of forcing
$\mathbb{P}$ which adjoins a graph $\Gamma$ of cardinality $\kappa$ such
that $G \simeq \Aut \Gamma$.
\end{enumerate}

Of course, there are many groups $G$ for which such a notion of forcing
$\mathbb{P}$ cannot possibly exist. For example, De Bruijn \cite{b} has shown  that 
the alternating group $\Alt ( \kappa^{+})$ cannot be embedded in $\Sym (\kappa)$.
Consequently, there is no cardinal-preserving notion of forcing 
$\mathbb{P}$ which adjoins a graph $\Gamma$ of cardinality $\kappa$ such
that $\Alt (\kappa^{+}) \simeq \Aut \Gamma$.

Our next definition singles out a combinatorial condition which is satisfied by
all those groups $G$ such that $G$ is embeddable in $\Sym (\kappa)$.
(See Proposition \ref{P:compatible}.) Conversely, in Theorem \ref{T:main2},
we will show that if a group $G$
satisfies this combinatorial condition, then there exists
a cardinal-preserving notion of forcing
$\mathbb{P}$ which adjoins a graph $\Gamma$ of cardinality $\kappa$ such
that $G \simeq \Aut \Gamma$.

\begin{Def} \label{D:compatible}
Let $\kappa$ be a regular uncountable cardinal such that $\kappa^{< \kappa} = \kappa$.
Then a group $G$ is said to satisfy the {\em $\kappa^{+}$-compatibility condition\/}
if it has the following property. Suppose that $H$ is a group such that
$\left| H \right| < \kappa$. Suppose that
$\left( f_{i} \mid i < \kappa^{+} \right)$ is a sequence of embeddings
$f_{i} : H \to G$; and let $H_{i} = f_{i} [H]$ for each
$i < \kappa^{+}$. Then there exist ordinals $i < j < \kappa^{+}$ and
a surjective homomorphism 
$\varphi : \langle H_{i}, H_{j} \rangle \to H_{i}$ such that
\begin{enumerate}
\item[(a)] $\varphi \circ f_{j} = f_{i}$; and
\item[(b)] $\varphi \res H_{i} = id_{H_{i}}$.
\end{enumerate}
\end{Def}

\begin{Example} \label{E:alt}
To get an understanding of Definition \ref{D:compatible},
it will probably be helpful to see an example of a group which
{\em fails\/} to satisfy the $\kappa^{+}$-compatibility condition. So
we will show that $\Alt (\kappa^{+})$ does not
satisfy the $\kappa^{+}$-compatibility condition.
Let $H = \Alt (4)$. For each $3 \leq i < \kappa^{+}$, let
$\Delta_{i} = \{ 0,1,2,i \}$ and let $f_{i} : \Alt (4) \to \Alt ( \Delta_{i})$ be an
isomorphism. If $3 \leq i < j < \kappa^{+}$, then
\[
\langle \Alt (\Delta_{i}) , \Alt (\Delta_{j}) \rangle =
\Alt (\Delta_{i} \cup \Delta_{j}) \simeq \Alt (5).
\]
Since $\Alt (5)$ is a simple group, there does not exist a surjective homomorphism
from $\langle \Alt (\Delta_{i}) , \Alt (\Delta_{j}) \rangle$ onto $\Alt ( \Delta_{i})$.
\end{Example}

\begin{Prop} \label{P:compatible}
Let $\kappa$ be a regular uncountable cardinal such that $\kappa^{< \kappa} = \kappa$,
and let $G \leqslant \Sym ( \kappa)$. Then $G$ satisfies the
$\kappa^{+}$-compatibility condition.
\end{Prop}

\begin{proof}
Let $H$ be a group such that $\left| H \right| < \kappa$, and let
$\left( f_{i} \mid i < \kappa^{+} \right)$ be a sequence of embeddings
$f_{i} : H \to G$. For each $i < \kappa^{+}$, let 
$H_{i} = f_{i} [H]$; and let $Z_{i}$ be a subset of $\kappa$ chosen so that
\begin{enumerate}
\item[(a)] $\left| Z_{i} \right| < \kappa$;
\item[(b)] $Z_{i}$ is $H_{i}$-invariant; and
\item[(c)] $g \res Z_{i} \neq id_{Z_{i}}$ for all $1 \neq g \in H_{i}$.
\end{enumerate}
After passing to a suitable subsequence if necessary, we can suppose that
the following conditions hold.
\begin{enumerate}
\item[(1)] There exists a fixed subset $Z$ such that $Z_{i} = Z$ for all
$i < \kappa^{+}$.
\item[(2)] For each $i < \kappa^{+}$, let
$r_{i} : H_{i} \to \Sym (Z)$ be the restriction mapping,
$g \mapsto g \res Z$; and let $\pi_{i} : H \to \Sym (Z)$ be the embedding
defined by $\pi_{i} = r_{i} \circ f_{i}$. Then $\pi_{i} = \pi_{j}$ for all
$i < j < \kappa^{+}$.
\end{enumerate}

Fix any pair of ordinals $i$, $j$ such that $i < j < \kappa^{+}$. Let
$\rho : \langle H_{i} , H_{j} \rangle \to \Sym (Z)$ be the restriction
mapping, $g \mapsto g \res Z$. Then
$\rho \left[\langle H_{i} , H_{j} \rangle \right] = 
r_{i} \left[ H_{i} \right]$, and so we can define a surjective homomorphism
$\varphi : \langle H_{i} , H_{j} \rangle \to H_{i}$ by
$\varphi = r_{i}^{-1} \circ \rho$. Clearly
$\varphi \res H_{i} = id_{H_{i}}$; and it is easily checked that
$\varphi \circ f_{j} = f_{i}$.
\end{proof}

\begin{Thm} \label{T:main2}
Let $\kappa$ be a regular uncountable cardinal such that $\kappa^{< \kappa} = \kappa$,
and let $G$ be a group which satisfies the $\kappa^{+}$-compatibility condition.
Then there exists a notion of forcing $\mathbb{P}$ such that
\begin{enumerate}
\item[(a)] $\mathbb{P}$ is $\kappa$-closed;
\item[(b)] $\mathbb{P}$ has the $\kappa^{+}$-$\,c.c.$; and
\item[(c)] $\underset{\mathbb{P}}{\Vdash} \text{ There exists a graph } \Gamma
\text{ of cardinality } \kappa \text{ such that } G \simeq \Aut \Gamma$.
\end{enumerate}
Furthermore, if $\left|G \right| = \lambda$, then
$\left| \mathbb{P} \right| = \max \{ \kappa , \lambda^{< \kappa} \}$.
\end{Thm}

Combining Proposition \ref{P:compatible} and Theorem \ref{T:main2},
we see that if $\kappa$ is an uncountable cardinal such that
$\kappa^{< \kappa} = \kappa$ and $G$ is an arbitrary subgroup
of $\Sym (\kappa)$, then there exists a cardinal-preserving notion
of forcing $\mathbb{P}$ and a graph $\Gamma \in V^{\mathbb{P}}$
such that $G \simeq \Aut \Gamma$. This result is {\em not\/} true
of arbitrary subgroups of $\Sym (\omega)$; for Solecki \cite{so}
has shown that no uncountable free abelian group is the 
automorphism group of a countable first-order structure.

Theorem \ref{T:main2} will be proved in Section \ref{S:closed}.
It is now easy to explain the main points of the proof of Theorem
\ref{T:best}. Assume that $V \vDash GCH$.
Let $\kappa$ be a regular uncountable cardinal, and let
$\lambda$ be a cardinal such that $\cf (\lambda) > \kappa$.
Let $\alpha$ be any ordinal such that $\alpha < \lambda^{+}$. 
In Section \ref{S:normaliser}, we will prove that
there exists a notion of forcing $\mathbb{Q}$ such that
\begin{enumerate}
\item[(1)] $\mathbb{Q}$ is $\kappa$-closed;
\item[(2)] $\mathbb{Q}$ has the $\kappa^{+}$-$\,c.c.$; 
\end{enumerate}
and such that the following statements are true in the generic
extension $V^{\mathbb{Q}}$.
\begin{enumerate}
\item[(a)] $2^{\kappa} = \lambda$.
\item[(b)] There exist groups $H < G < \Sym( \kappa)$ such that 
$|H| = \kappa$ and the normaliser tower of $H$ in $G$ terminates
in exactly $\alpha$ steps.
\end{enumerate}
By Proposition \ref{P:compatible}, $G$
satisfies the $\kappa^{+}$-compatibility condition. Hence we
can use Theorem \ref{T:main2} to generically adjoin a graph
$\Gamma$ of cardinality $\kappa$ such that $G \simeq 
\Aut \Gamma$. A moment's thought shows that the normaliser tower
of $H$ in $G$ is an absolute notion. Thus
Lemma \ref{L:reduction} yields a centreless
group $T$ of cardinality $\kappa$ such that $\tau (T) = \alpha$.
The case when $\kappa$ is a singular cardinal requires a
little more work, and will be dealt with in Section 
\ref{S:increasing}.
The remainder of this section will be devoted to another two easy
applications of Theorem \ref{T:main2}.

\begin{Application} \label{A:free}
A well-known open problem asks whether there exists a countable
structure $\mathcal{M}$ such that $\Aut \mathcal{M}$ is the
free group on $2^{\omega}$ generators. Using Theorem
\ref{T:main2}, it is easy to establish the consistency of the
existence of a structure $\mathcal{N}$ of cardinality $\omega_{1}$
such that $\Aut \mathcal{N}$ is the
free group on $2^{\omega_{1}}$ generators. It is not known whether
the existence of such a structure can be proved in $ZFC$.

\begin{Thm} \label{T:free}
Let $V$ be a transitive model of $ZFC$ and let $\kappa$, $\lambda$,
$\theta$ be cardinals such that
$\kappa^{< \kappa} = \kappa < \lambda \leq
\theta = \theta^{\kappa}$. Then there exists a notion of forcing
$\mathbb{P}$, which preserves cofinalities and cardinalities,
such that the following statements are true in $V^{\mathbb{P}}$.
\begin{enumerate}
\item[(a)] $2^{\kappa} = \theta$; and
\item[(b)] there exists a graph $\Gamma$ of cardinality $\kappa$
such that $\Aut \Gamma$ is the free group on $\lambda$
generators.
\end{enumerate}
\end{Thm}

\begin{proof}
After performing a preliminary forcing if necessary, we can also suppose
that $2^{\kappa} = \theta$. Let $F$ be the free group on $\lambda$
generators. Then it is enough to show that $F$ satisfies the
$\kappa^{+}$-compatibility condition. Let $H$ be a (necessarily free) group such
that $\left| H \right| < \kappa$, and let 
$\left( f_{i} \mid i < \kappa^{+} \right)$ be a sequence of embeddings
$f_{i} : H \to F$. For each $i < \kappa^{+}$, let 
$H_{i} = f_{i} [H]$. Then
$\langle H_{i} \mid i < \kappa^{+} \rangle$ is a free group of cardinality at most
$\kappa^{+}$. By Theorem 5.1 \cite{s}, there exists an embedding of
$\langle H_{i} \mid i < \kappa^{+} \rangle$ into $\Sym (\kappa)$. So Proposition
\ref{P:compatible} yields the existence of ordinals $i < j < \kappa^{+}$ and a 
surjective homomorphism 
$\varphi : \langle H_{i}, H_{j} \rangle \to H_{i}$ such that
$\varphi \circ f_{j} = f_{i}$ and $\varphi \res H_{i} = id_{H_{i}}$.
\end{proof}
\end{Application}

\begin{Application} \label{A:kueker}
Theorem 1.6 \cite{t2} says that if $G$ is a finitely generated centreless
group, then the automorphism tower of $G$ terminates in countably
many steps. It is conceivable that a more general result holds; namely,
that the automorphism tower of $G$ terminates in countably many
steps, whenever $G$ is a countable centreless group such that
$\Aut G$ is also countable. To see why this might be true, let $G$ be
such a group. Then, by Kueker \cite{ku}, there exists a finite subset
$F \subseteq G$ such that each automorphism $\pi \in \Aut G$ is
uniquely determined by its restriction $\pi \res F$. In terms of the
automorphism tower of $G$, this says that there is a finite subset
$F \subseteq G$ such that $C_{G_{1}}(F) = 1$. Suppose that the
``rigidity'' of $F$ within $G = G_{0}$ is propagated along the
automorphism tower of $G$; ie. that 
$C_{G_{\alpha}}(F) =1$ for all ordinals $\alpha < \omega_{1}$. Then the proof 
of Theorem 1.6 \cite{t2} shows that the automorphism tower of
$G$ terminates in countably many steps.

\begin{Question} \label{Q:kueker}
Let $G$ be a centreless group such that 
$\left| \Aut G \right| = \omega$. Does there exist a finite subset
$F \subseteq G$ such that 
$C_{G_{\alpha}}(F) =1$ for all ordinals $\alpha < \omega_{1}$?
\end{Question}

If $\left| G_{\alpha} \right| = \omega$ for all $\alpha < \omega_{1}$,
then Fodor's Lemma implies that there exists an
ordinal $\beta < \omega_{1}$ and a finite subset $F_{\beta}$ of
$G_{\beta}$ such that $C_{G_{\alpha}}(F_{\beta}) = 1$ for all
$\beta \leq \alpha < \omega_{1}$; and so $\tau(G) < \omega_{1}$.
This observation suggests the following weak form of 
Question \ref{Q:kueker}, which is also open.

\begin{Question} \label{Q:kueker2}
Does there exist a centreless group $G$ such that
$\left| \Aut G \right| = \omega$ and
$\left| \Aut ( \Aut G) \right| = 2^{\omega}$?
\end{Question}

Of course, a positive answer to Question \ref{Q:kueker} implies a negative
answer to Question \ref{Q:kueker2}. Using Theorem \ref{T:main2}, it is
easy to establish the consistency of the existence of a centreless group
$G$ of cardinality $\omega_{1}$ such that
$\left| \Aut G \right| = \omega_{1}$ and
$\left| \Aut ( \Aut G) \right| = 2^{\omega_{1}}$. Once again, it is not
known whether the existence of such a group can be proved in $ZFC$.

\begin{Thm} \label{T:kueker}
Let $\kappa$ be a regular uncountable cardinal such that
$\kappa^{< \kappa} = \kappa$. Then it is consistent that there exists
a centreless group $G$ of cardinality $\kappa$ such that
$\left| \Aut G \right| = \kappa$ and
$\left| \Aut ( \Aut G) \right| = 2^{\kappa}$.
\end{Thm}

\begin{proof}
Let $V$ be the ground model. For each $\alpha$, $\xi < \kappa$,
let $Z^{\alpha}_{\xi} = \langle z^{\alpha}_{\xi} \rangle$ be an infinite
cyclic group. For each $\xi < \kappa$, let
$A_{\xi} = \underset{\alpha < \kappa}{\bigoplus}Z^{\alpha}_{\xi}$;
and let $B = \underset{\xi < \kappa}{\bigoplus}A_{\xi}$. Define
an action of $\Sym (\kappa)$ on $B$ by
$\pi z^{\alpha}_{\xi} \pi^{-1} = z^{\alpha}_{ \pi (\xi)}$ for all
$\alpha$, $\xi < \kappa$; and let
$W = B \rtimes \Sym (\kappa)$ be the corresponding semidirect
product. Let $H = \underset{\xi < \kappa}{\bigoplus}
Z^{\xi}_{\xi}$. Then the members of the normaliser tower of
$H$ in $W$ are
\begin{enumerate}
\item[(a)] $N_{0}(H) = H$;
\item[(b)] $N_{1}(H) = B$;
\item[(c)] $N_{2}(H) = W$.
\end{enumerate}
Clearly $W$ is embeddable in $\Sym (\kappa)$; and so $W$ satisfies the
$\kappa^{+}$-compatibility condition. Let $\mathbb{P}$ be the notion of 
forcing, given by Theorem \ref{T:main2}, which adjoins a graph 
$\Gamma$ of cardinality $\kappa$ such that $W \simeq \Aut \Gamma$.
Let $K_{\Gamma} \in V^{\mathbb{P}}$ be the corresponding field,
which is given by Lemma \ref{L:field}. Then
$G = PGL( 2, K_{\Gamma}) \rtimes H$ is a group such that
\[
\left| \Aut G \right| = \left| PGL( 2, K_{\Gamma}) \rtimes B \right|
= \kappa
\]
and
\[
\left| \Aut ( \Aut G) \right| = \left| PGL( 2, K_{\Gamma}) \rtimes W \right|
= 2^{\kappa}.
\]
\end{proof}
\end{Application}

Our set-theoretic notation mainly follows that 
of Kunen \cite{k}. Thus if $\mathbb{P}$
is a notion of forcing and $p$,$q \in \mathbb{P}$, then $q \leq p$
means that $q$ is a strengthening of $p$. We say that $\mathbb{P}$
is $\kappa$-closed if for every $\lambda < \kappa$, every 
descending sequence of elements of $\mathbb{P}$
\[
p_{0} \geq p_{1} \geq \dots \geq p_{\xi} \geq \dots , \qquad  \xi < \lambda ,
\]
has a lower bound in $\mathbb{P}$. If $V$ is the ground model,
then we will denote the generic extension by $V^{\mathbb{P}}$ if
we do not wish to specify a particular generic filter 
$H \subseteq \mathbb{P}$. If we want to emphasize that the term $t$
is to be interpreted in the model $M$ of $ZFC$, then we write
$t^{M}$; for example, $\Sym (\lambda)^{M}$.

Our group-theoretic notation is standard. For example, 
the (restricted) wreath product of $A$ by $C$ is denoted by
$A \Wr C$; and the direct sum of the groups
$H_{\xi}$, $\xi < \lambda$, is denoted by
$\underset{\xi < \lambda}{\bigoplus} H_{\xi}$. If 
$\pi \in \Sym ( \kappa)$, then
$\supp (\pi) = \{ \alpha \in \kappa \mid \pi (\alpha) \neq \alpha \}$;
and if $\lambda$ is an infinite cardinal such that 
$\lambda \leq \kappa$, then
$\Sym_{\lambda}(\kappa) = \{ \pi \in \Sym (\kappa) \mid
\left| \supp (\pi) \right| < \lambda \}$.

\section{Realising normaliser towers within infinite symmetric groups} \label{S:normaliser}

Let $\kappa$ be a regular uncountable cardinal such that 
$\kappa^{< \kappa} = \kappa$. In this section, we will study the problem
of realising long normaliser towers within $\Sym (\kappa)$. In particular,
we will prove that if 
$\alpha$ is any ordinal, then there exists a generic extension
$V^{\mathbb{Q}}$ such that a normaliser tower of height $\alpha$
can be realised in $\Sym (\kappa)^{V^{\mathbb{Q}}}$.

First for each ordinal $\alpha$, we will construct a pair of
groups, $H < G$, such that the normaliser tower of $H$ in $G$
terminates in exactly $\alpha$ steps.

\begin{Def} \label{D:wreath}
The ascending chain of groups
\[
W_{0} \leqslant W_{1} \leqslant \dots \leqslant
W_{\alpha} \leqslant W_{\alpha +1} \leqslant \dots
\]
is defined inductively as follows.
\begin{enumerate}
\item[(a)] $W_{0} = C_{2}$, the cyclic group of order 2.
\item[(b)] Suppose that $\alpha = \beta +1$. Then
\[
W_{\beta} = W_{\beta} \oplus 1 \leqslant
\left[ W_{\beta} \oplus W_{\beta}^{*} \right]
\rtimes \langle \sigma_{\beta +1} \rangle =
W_{\beta + 1}.
\]
Here $W_{\beta}^{*}$ is an isomorphic copy of
$W_{\beta}$; and $\sigma_{\beta +1}$ is an element
of order 2 which interchanges the factors
$W_{\beta} \oplus 1$ and $1 \oplus W_{\beta}^{*}$
of the direct sum
$ W_{\beta} \oplus W_{\beta}^{*}$ via conjugation.
Thus $W_{\beta +1}$ is isomorphic to the wreath product
$W_{\beta} \Wr C_{2}$.
\item[(c)] If $\alpha$ is a limit ordinal, then
$W_{\alpha} = \underset{\beta < \alpha}{\bigcup}
W_{\beta}$.
\end{enumerate}
\end{Def}

\begin{Lem} \label{L:wreath}
$\left| W_{\alpha} \right| \leq
\max \{ |\alpha|, \omega \}$ for all ordinals
$\alpha$.
\end{Lem}

\begin{proof}
This follows by an easy induction on $\alpha$.
\end{proof}

\begin{Lem} \label{L:wreath1}
\begin{enumerate}
\item[(a)] If $1 \leq n < \omega$, then the normaliser tower of $W_{0}$ 
in $W_{n}$ terminates in exactly $n+1$ steps.
\item[(b)] If $\alpha \geq \omega$, then the normaliser tower of $W_{0}$ 
in $W_{\alpha}$ terminates in exactly $\alpha$ steps.
\end{enumerate}
\end{Lem}

\begin{proof}
(a) It is easily checked that 
\[
N_{1}(W_{0},W_{n}) = W_{0} \oplus W_{0}^{*} \oplus W_{1}^{*}
\oplus \dots \oplus W_{n-1}^{*}
\]
and that
\[
N_{2}(W_{0},W_{n}) = W_{1}  \oplus W_{1}^{*}
\oplus \dots \oplus W_{n-1}^{*};
\]
and that, in general, for each $0 \leq \ell \leq n-1$,
\[
N_{\ell +1}(W_{0},W_{n}) 
= W_{\ell} \oplus \underset{\ell \leq m < n}{\bigoplus} W_{m}^{*}.
\]

(b) For example, consider the case when $\alpha > \omega$. Then
for each $\ell \in \omega$,
\[
N_{\ell +1}(W_{0},W_{\alpha}) 
= W_{\ell} \oplus \underset{\ell \leq \beta < \alpha}{\bigoplus} 
W_{\beta}^{*};
\]
and for each $\gamma$ such that $\omega \leq \gamma < \alpha$,
\[
N_{\gamma}(W_{0},W_{\alpha}) 
= W_{\gamma} \oplus \underset{\gamma \leq \beta < \alpha}{\bigoplus} 
W_{\beta}^{*}.
\]
\end{proof}

\begin{Remark} \label{R:wreath}
Unfortunately, the group $W_{\kappa^{+}}$ does not satisfy
the $\kappa^{+}$-compatibility condition. To see this, let
\[
H = C_{2} \Wr C_{2} = \left[ \langle a \rangle \oplus
\langle b \rangle \right] \rtimes \langle c \rangle ;
\]
and for each limit ordinal $i < \kappa^{+}$, let
$f_{i} : H \to W_{\kappa^{+}}$ be the embedding
such that $f_{i}(a) = \sigma_{i+1}$ and
$f_{i}(c) = \sigma_{i+2}$. (Here we are using the notation
which was introduced in Definition \ref{D:wreath}.) Let
$i$, $j$ be any limit ordinals such that 
$i < j < \kappa^{+}$. Then
\[
\langle H_{i}, H_{j} \rangle \simeq
\left( \left( C_{2} \Wr C_{2} \right) \Wr C_{2} \right) 
\Wr C_{2}.
\]
Suppose that there exists a a surjective homomorphism 
$\varphi : \langle H_{i}, H_{j} \rangle \to H_{i}$ such that
\begin{enumerate}
\item[(a)] $\varphi \circ f_{j} = f_{i}$; and
\item[(b)] $\varphi \res H_{i} = id_{H_{i}}$.
\end{enumerate}
Then 
$\varphi(\sigma_{j+1}) = \varphi(\sigma_{i+1}) = \sigma_{i+1}$
and
$\varphi(\sigma_{j+2}) = \varphi(\sigma_{i+2}) = \sigma_{i+2}$.
Consider the element 
$x = \sigma_{j+1} \sigma_{j+2} \sigma_{j+1} \sigma_{j+2}
\in H_{j}$. Then it is easily checked that
\begin{enumerate}
\item[(1)] $x$ lies in the centre of $H_{j}$; and
\item[(2)] $\sigma_{j+1} y \sigma_{j+1}^{-1} =
x y x^{-1}$ for all $y \in H_{i}$.
\end{enumerate}
Thus $z = \varphi(x)$ lies in the centre of $H_{i}$. Since
$\varphi(\sigma_{j+1}) = \sigma_{i+1}$, we find that
\[
\sigma_{i+1} y \sigma_{i+1}^{-1} = z y z^{-1} = y
\]
for all $y \in H_{i}$. But this contradicts the fact that
$\sigma_{i+1}$ is a noncentral element of $H_{i}$.
\end{Remark}

Thus if $\alpha \geq \kappa^{+}$, then
$W_{\alpha}$ is not embeddable in
$\Sym (\kappa)$. However, the above argument does not rule out
the possibility that $W_{\alpha}$ is embeddable in the quotient
group $\Sym (\kappa) / \Sym_{\kappa}(\kappa)$; and this is
enough for our purposes.

\begin{Lem} \label{L:enough}
Let $\kappa$ be an infinite cardinal such that
$\kappa^{< \kappa} = \kappa$. Suppose that $\gamma$ is an
ordinal, and that there exists an embedding
\[
f : W_{\gamma} \to \Sym (\kappa)/ \Sym_{\kappa}(\kappa).
\]
Then for each ordinal $\alpha \leq \gamma$, there exist groups
$H \leqslant G \leqslant \Sym(\kappa)$ such that
$|H| = \kappa$ and the normaliser tower of $H$ in $G$
terminates in exactly $\alpha$ steps.
\end{Lem}

\begin{proof}
For each $\alpha \leq \gamma$, let $H^{\alpha}$ be the subgroup
of $\Sym(\kappa)$ such that $f \left[ W_{\alpha} \right] =
H^{\alpha} / \Sym_{\kappa}(\kappa)$. Since
$| \Sym_{\kappa}(\kappa)| = \kappa^{< \kappa} = \kappa$, it 
follows that $|H^{0}|= \kappa$.

\begin{Claim} \label{C:enough}
For each ordinal $\beta$,
\[
f \left[ N_{\beta}(W_{0},W_{\alpha}) \right] =
N_{\beta}(H^{0},H^{\alpha}) / \Sym_{\kappa}(\kappa).
\]
\end{Claim}

\begin{proof}
We will argue by induction on $\beta$. The result is clear when
$\beta = 0$, and no difficulties arise when $\beta$ is a limit
ordinal. Suppose that $\beta = \xi +1$ and that the result
holds for $\xi$. Let $R = N_{\xi}(H^{0},H^{\alpha})$; and for
each subgroup $K$ such that
$\Sym_{\kappa}(\kappa) \leqslant K \leqslant H^{\alpha}$, let
$\overline{K} = K/ \Sym_{\kappa}(\kappa)$. Then
\[
f \left[ N_{\xi +1}(W_{0},W_{\alpha}) \right] =
N_{\overline{H^{\alpha}}}(\overline{R});
\]
and so we must show that
\[
N_{\overline{H^{\alpha}}}(\overline{R}) =
\overline{N_{H^{\alpha}}(R)}.
\]
But this is an immediate 
consequence of the Correspondence Theorem for subgroups of
quotient groups, together with the observation that
the normaliser of any subgroup $L$ is the largest subgroup
$M$ such that $L \trianglelefteq M$.
\end{proof}

It is now easy to complete the proof of Lemma \ref{L:enough}.
Applying Lemma \ref{L:wreath1} and Claim \ref{C:enough}, we see that
if $\alpha \geq \omega$, then the normaliser tower of $H^{0}$ in
$H^{\alpha}$ terminates in exactly $\alpha$ steps; and that if
$2 \leq \alpha = n < \omega$, then the normaliser tower of
$H^{0}$ in $H^{n-1}$ terminates in exactly $n$ steps. This
just leaves the cases when $\alpha = 0, 1$. When $\alpha =0$,
we can take $H = G = \Alt (\kappa)$; and when $\alpha =1$, we
can take $H = \Alt (\kappa)$ and $G = \Sym(\kappa)$.
\end{proof}

The next result implies that if
$\omega < \kappa^{\kappa} = \kappa$ and $W$ is {\em any\/}
group, then there exists a cardinal-preserving notion of
forcing $\mathbb{P}$ such that in $V^{\mathbb{P}}$, the
group $W$ is embeddable in 
$\Sym (\kappa)/ \Sym_{\kappa}(\kappa)$. (If 
$|W|^{\kappa} > |W|$, then just embed $W$ in a larger group
$L$ such that $|L|^{\kappa} = |L|$; and then apply
Lemma \ref{L:just} to $L$.)

\begin{Lem} \label{L:just}
Let $V$ be a transitive model of $ZFC$ and let $\kappa$, $\lambda$
be cardinals such that 
$\omega < \kappa^{<\kappa}=\kappa < \lambda = \lambda^{\kappa}$.
Let $W$ be any group of cardinality $\lambda$.
Then there exists a notion of forcing $\mathbb{Q}$ such that
\begin{enumerate}
\item[(1)] $\mathbb{Q}$ is $\kappa$-closed.
\item[(2)] $\mathbb{Q}$ has the $\kappa^{+}$-$\,c.c.$; 
\end{enumerate}
and such that the following statements are true in the generic
extension $V^{\mathbb{Q}}$.
\begin{enumerate}
\item[(a)] $2^{\kappa} = \lambda$.
\item[(b)] There exists an isomorphic embedding
\[
f : W \to 
\left( \Sym (\kappa)/ \Sym_{\kappa}(\kappa) \right)^{V^{\mathbb{Q}}}.
\]
\end{enumerate}
\end{Lem}

\begin{proof}
Let $\Omega = \bigcup_{\alpha < \kappa}\{\alpha \} \times \alpha$.
We will work with the symmetric group $\Sym (\Omega)$ rather
than with $\Sym (\kappa)$. Let $\mathbb{Q}$ be the notion of
forcing consisting of the conditions
\[
p = ( \delta_{p}, H_{p}, E_{p} )
\]
such that the following hold.
\begin{enumerate}
\item[(a)] $\omega \leq \delta_{p} < \kappa$.
\item[(b)] $H_{p}$ is a subgroup of $W$
such that $|H_{p}| \leq |\delta_{p}|$.
\item[(c)] $E_{p}$ is a function which assigns a permutation
$e^{p}_{\pi, \xi} \in \Sym ( \{ \xi \} \times \xi )$ to each
pair $(\pi, \xi) \in H_{p} \times \delta_{p}$.
\end{enumerate}
We set $q = ( \delta_{q}, H_{q}, E_{q} ) \leq
p = ( \delta_{p}, H_{p}, E_{p} )$ if and only if
\begin{enumerate}
\item[(1)] $\delta_{p} \leq \delta_{q}$;
\item[(2)] $H_{p} \leqslant H_{q}$;
\item[(3)] $E_{p} \subseteq E_{q}$; and
\item[(4)] if $\delta_{p} \leq \xi < \delta_{q}$, then the
restriction to $H_{p}$ of the function,
$\pi \mapsto e^{q}_{\pi, \xi}$, is an isomorphic embedding
of $H_{p}$ into $\Sym (\{\xi\} \times \xi)$.
\end{enumerate}

\begin{Claim} \label{C:just1}
$\mathbb{Q}$ is $\kappa$-closed.
\end{Claim}

\begin{proof}[Proof of Claim \ref{C:just1}]
This is clear.
\end{proof}

\begin{Claim} \label{C:just2}
$\mathbb{Q}$ has the $\kappa^{+}$-$\,c.c.$.
\end{Claim}

\begin{proof}[Proof of Claim \ref{C:just2}]
Suppose that 
$p_{i} = ( \delta_{p_{i}}, H_{p_{i}}, E_{p_{i}} ) \in \mathbb{Q}$
for $i < \kappa^{+}$. Using the $\Delta$-System Lemma and the
assumption that $\kappa^{<\kappa}= \kappa$, after passing to
a suitable subsequence if necessary, we can suppose that the
following conditions are satisfied.
\begin{enumerate}
\item[(i)] There exists a fixed ordinal $\delta$ such that
$\delta_{p_{i}} = \delta$ for all $i < \kappa^{+}$.
\item[(ii)] There exists a fixed subgroup $H$ such that
$H_{p_{i}} \cap H_{p_{j}} = H$ for all $i < j < \kappa^{+}$.
\item[(iii)] There exists a fixed function $E$ such that
$E_{p_{i}} \res H = E$ for all $i < \kappa^{+}$.
\end{enumerate}
Now fix any two ordinals $i < j < \kappa^{+}$. Let
$H^{+} = \langle H_{p_{i}}, H_{p_{j}} \rangle$ be the subgroup
generated by $H_{p_{i}} \cup H_{p_{j}}$; and let $E^{+}$ be any
extension of $E_{p_{i}} \cup E_{p_{j}}$ which satisfies
condition (c). Then $q = (\delta, H^{+}, E^{+} )$ is a common
lower bound of $p_{i}$ and $p_{j}$.
\end{proof}

\begin{Claim} \label{C:just3}
For each $\alpha < \kappa$, the set
$C_{\alpha} = \{ p \in \mathbb{Q} \mid \delta_{p} \geq \alpha \}$
is dense in $\mathbb{Q}$.
\end{Claim}

\begin{proof}[Proof of Claim \ref{C:just3}]
Let $p = ( \delta_{p}, H_{p}, E_{p} ) \in \mathbb{Q}$.
Then we can suppose that $\delta_{p} < \alpha$. We can
define an isomorphic embedding
$\varphi : H_{p} \to \Sym (H_{p})$ by setting
$\varphi (h)(x) = hx$ for all $x \in H_{p}$. Since
$\left|H_{p} \right| \leq \left| \delta_{p} \right|$,
it follows that there exists an isomorphic embedding
$\varphi_{\xi} : H_{p} \to \Sym (\{\xi \} \times \xi)$
for each $\delta_{p} \leq \xi < \alpha$. Hence there
exists a condition
$q = ( \delta_{q}, H_{q}, E_{q} ) \leq p$ such that
$H_{q} = H_{p}$ and $\delta_{q} = \alpha$.
\end{proof}

\begin{Claim} \label{C:just4}
For each $\pi \in W$, the set
$D_{\pi} = \{ p \in \mathbb{Q} \mid \pi \in H_{p} \}$
is dense in $\mathbb{Q}$.
\end{Claim}

\begin{proof}[Proof of Claim \ref{C:just4}]
Let $p = ( \delta_{p}, H_{p}, E_{p} ) \in \mathbb{Q}$.
Let $H^{+} = \langle H_{p}, \pi \rangle$ be the 
subgroup generated by $H_{p} \cup \{ \pi \}$, and let
$E^{+}$ be any extension of $E_{p}$ to $H^{+}$ which
satisfies condition (c). Then
$q = ( \delta_{q}, H_{q}, E_{q} ) \leq p$.
\end{proof}

Let $F \subseteq \mathbb{Q}$ be a generic filter, and let
$V^{\mathbb{Q}} = V[F]$ be the corresponding generic
extension. Working within $V^{\mathbb{Q}}$, for each
$\pi \in W$, let
\[
e( \pi) = \bigcup \{ e^{p}_{\pi, \xi } \mid
\text{ There exists } p \in F \text{ such that }
\pi \in H_{p} \text{ and } \xi < \delta_{p} \}.
\]
Then $e( \pi) \in \Sym (\Omega)$. Let
$\Sym_{\kappa}(\Omega) = \{ \psi \in \Sym (\Omega) \mid
|\supp (\psi)| < \kappa \}$; and define the function
\[
f: W \to \Sym (\Omega)/ \Sym_{\kappa}(\Omega)
\]
by $f(\pi) = e(\pi) \Sym_{\kappa}(\Omega)$. Then it is enough to 
show that $f$ is an isomorphic embedding.

\begin{Claim} \label{C:just5}
If $1 \neq \pi \in W$, then $f(\pi) \neq 1$.
\end{Claim}

\begin{proof}[Proof of Claim \ref{C:just5}]
Choose a condition 
$p = ( \delta_{p}, H_{p}, E_{p} ) \in F$ such that $\pi \in H_{p}$. 
If $\xi$ is any ordinal such that $\delta_{p} \leq \xi < \kappa$,
then $e(\pi) \res \{ \xi \} \times \xi \neq id_{\{ \xi \} \times \xi}$.
Hence $e(\pi) \notin \Sym_{\kappa}(\Omega)$.
\end{proof}

\begin{Claim} \label{C:just6}
$f$ is a group homomorphism.
\end{Claim}

\begin{proof}[Proof of Claim \ref{C:just6}]
Let $\pi_{1}$, $\pi_{2} \in W$.
Let $p = ( \delta_{p}, H_{p}, E_{p} ) \in F$ be
a condition such that 
$\pi_{1}$, $\pi_{2} \in H_{p}$. Let $\xi$ be any ordinal
such that $\delta_{p} \leq \xi < \kappa$; and let
$q \in F$ be a condition such that $q \leq p$ and
$\xi < \delta_{q}$. Then
\[
e^{q}_{\pi_{1},\xi} \circ e^{q}_{\pi_{2},\xi} =
e^{q}_{\pi_{1} \circ \pi_{2},\xi}.
\]
It follows that
\[
e(\pi_{1}) \Sym_{\kappa}(\Omega) \circ
e(\pi_{2}) \Sym_{\kappa}(\Omega) =
e(\pi_{1} \circ \pi_{2}) \Sym_{\kappa}(\Omega).
\]
\end{proof}

Finally it is easily checked that $|\mathbb{Q}| = \lambda$;
and it follows that 
$V^{\mathbb{Q}} \vDash 2^{\kappa} = \lambda$.
This completes the proof of Lemma \ref{L:just}.
\end{proof}

Summing up our work in this section, we obtain the following 
result.

\begin{Thm} \label{T:just}
Let $V$ be a transitive model of $ZFC$ and let $\kappa$, $\lambda$
be cardinals such that 
$\omega < \kappa^{<\kappa}=\kappa < \lambda = \lambda^{\kappa}$.
Let $\alpha$ be any ordinal such that $\alpha < \lambda^{+}$.
Then there exists a notion of forcing $\mathbb{Q}$ such that
\begin{enumerate}
\item[(1)] $\mathbb{Q}$ is $\kappa$-closed;
\item[(2)] $\mathbb{Q}$ has the $\kappa^{+}$-$\,c.c.$; 
\end{enumerate}
and such that the following statements are true in the generic
extension $V^{\mathbb{Q}}$.
\begin{enumerate}
\item[(a)] $2^{\kappa} = \lambda$.
\item[(b)] There exist groups $H < G < \Sym( \kappa)$ such that 
$|H| = \kappa$ and the normaliser tower of $H$ in $G$ terminates
in exactly $\alpha$ steps.
\end{enumerate}
\end{Thm}

\begin{proof}
Let $\gamma$ be an ordinal such that
$\max \{ \alpha, \lambda \} \leq \gamma < \lambda^{+}$. Then
$|W_{\gamma}| = \lambda$. Let $\mathbb{Q}$ be the notion of forcing
obtained by applying Lemma \ref{L:just} to
$W = W_{\gamma}$. By Lemma \ref{L:enough}, $\mathbb{Q}$ satisfies
our requirements.
\end{proof}

\section{Closed groups of uncountable degree} \label{S:closed}
In this section, we will prove Theorem \ref{T:main2}.
Let $\kappa$ be a regular uncountable cardinal such
that $\kappa^{<\kappa} = \kappa$, and let $G$ be a group which
satisfies the $\kappa^{+}$-compatibility condition. Let $L$ be a
first-order language consisting of $\kappa$ binary relation symbols.
The following notion of forcing $\mathbb{P}$ is designed to adjoin
a structure $\mathcal{M}$ of cardinality $\kappa$ for the language
$L$ such that $G \simeq \Aut \mathcal{M}$. This is sufficient; for then
we can use one of the standard procedures to code $\mathcal{M}$
into a graph $\Gamma$ of cardinality $\kappa$ such that
$\Aut \Gamma \simeq \Aut \mathcal{M}$. ( Cf. Section 5.5 of
Hodges \cite{h}.)

\begin{Def} \label{D:atomic}
Suppose that $L_{0} \subseteq L$ and that $\mathcal{N}$ is a
structure for the language $L_{0}$. Then a
{\em restricted atomic type\/} in the free variable $v$ for the
language $L_{0}$ using parameters from $\mathcal{N}$ is a
set $t$ of formulas of the form $R(v,a)$, where $R \in L_{0}$ 
and $a \in \mathcal{N}$. 
An element $c \in \mathcal{N}$ is said to {\em realise\/} $t$ if
$\mathcal{N} \vDash \varphi[c]$ for every formula
$\varphi(v) \in t$. If no element of $\mathcal{N}$ realises
$t$, then $t$ is said to be {\em omitted\/} in $\mathcal{N}$. Notice
that if $t$ is omitted in $\mathcal{N}$, then $t$ is not the
{\em trivial\/} restricted atomic type $\emptyset$.
\end{Def}

\begin{Def} \label{D:poset}
Let $\mathbb{P}$ be the notion of forcing consisting of the
conditions
\[
p = \left( H, \pi, \mathcal{N}, T \right)
\]
such that the following hold.
\begin{enumerate}
\item[(a)] $H$ is a subgroup of $G$ such that
$\left| H \right| < \kappa$.
\item[(b)] There exists an ordinal $0 < \delta < \kappa$
and a subset $L(\mathcal{N}) \in [L]^{< \kappa}$ 
such that $\mathcal{N}$ is a
structure with universe $\delta$ for the 
language $L(\mathcal{N})$.
\item[(c)] $\pi: H \to \Aut \mathcal{N}$ is a group homomorphism.
\item[(d)] $T$ is a set of restricted atomic types 
in the free variable $v$ for the language
$L( \mathcal{N})$ using parameters from $\mathcal{N}$. Furthermore,
$\left| T \right| < \kappa$; and each $t \in T$ is omitted in $\mathcal{N}$.
\end{enumerate}
We set $\left( H_{2}, \pi_{2}, \mathcal{N}_{2}, T_{2} \right) \leq
\left( H_{1}, \pi_{1}, \mathcal{N}_{1}, T_{1} \right)$ if and only if
\begin{enumerate}
\item[(1)] $H_{1} \leqslant H_{2}$.
\item[(2)] $\mathcal{N}_{1}$ is a substructure of $\mathcal{N}_{2}$.
\item[(3)] For all $h \in H_{1}$ and $\alpha \in \mathcal{N}_{1}$,
$\pi_{2}(h)(\alpha) = \pi_{1}(h)(\alpha)$.
\item[(4)] $T_{1} \subseteq T_{2}$.
\end{enumerate}
\end{Def}

It is clear that the components $\left(H, \pi, \mathcal{N} \right)$
in each condition $p \in \mathbb{P}$ are designed to generically
adjoin a structure $\mathcal{M}$ of cardinality $\kappa$ for the
language $L$, together with an embedding $\pi^{*}$ of $G$ into
$\Aut \mathcal{M}$. The set $T$ of restricted atomic types is needed to 
kill off potential extra automorphisms
$g \in \Aut \mathcal{M} \smallsetminus \pi^{*}[G]$, and thus ensure
that $\pi^{*}$ is surjective.

\begin{Lem} \label{L:poset1}
$\mathbb{P}$ is $\kappa$-closed.
\end{Lem}

\begin{proof}
This is clear.
\end{proof}

\begin{Lem} \label{L:poset2}
$\mathbb{P}$ has the $\kappa^{+}$-$\, c.c.$
\end{Lem}

\begin{proof}
Suppose that
$ p_{i} = \left( H_{i}, \pi_{i}, \mathcal{N}_{i}, T_{i} \right) 
\in \mathbb{P}$ for $i < \kappa^{+}$. After passing to a suitable 
subsequence if necessary, we can suppose that the following 
conditions hold.
\begin{enumerate}
\item[(1)] There exists a fixed structure $\mathcal{N}$ such that
$\mathcal{N}_{i} = \mathcal{N}$ for all $i < \kappa^{+}$.
\item[(2)] There exists a fixed set of restricted atomic types $T$ such that
$T_{i} = T$ for all $i < \kappa^{+}$.
\item[(3)] There is a fixed group $H$ such that for each
$i < \kappa^{+}$, there exists an isomorphism
$f_{i} : H \simeq H_{i}$.
\item[(4)] For each $i < \kappa^{+}$, let
$\psi_{i} : H \to \Aut \mathcal{N}$ be the embedding
defined by $\psi_{i} = \pi_{i} \circ f_{i}$. Then
$\psi_{i} = \psi_{j}$ for all $i < j < \kappa^{+}$.
\end{enumerate}
Since $G$ satisfies the $\kappa^{+}$-compatibility
condition, there exist ordinals $i < j < \kappa^{+}$ and
a surjective homomorphism 
$\varphi : \langle H_{i}, H_{j} \rangle \to H_{i}$ such that
\begin{enumerate}
\item[(a)] $\varphi \circ f_{j} = f_{i}$; and
\item[(b)] $\varphi \res H_{i} = id_{H_{i}}$.
\end{enumerate}
Let $\langle H_{i},H_{j} \rangle \to \Aut \mathcal{N}$
be the homomorphism defined by
$\pi = \pi_{i} \circ \varphi$. Clearly $\pi_{i} \subseteq \pi$.
Note that if $x \in H_{j}$, then
\begin{align*}
\pi_{i} \circ \varphi (x) &= 
\pi_{i} \circ \left( \varphi \circ f_{j} \right) \circ f_{j}^{-1}(x) \\
&= \pi_{i} \circ f_{i} \circ f_{j}^{-1}(x) \\
&= \pi_{j} (x) .
\end{align*}
Thus we also have that $\pi_{j} \subseteq \pi$. Consequently, we can
define a condition $p \leq p_{i}$, $p_{j}$ by
\[
p = \left( \langle H_{i},H_{j} \rangle, \pi, \mathcal{N}, T \right) .
\]
\end{proof}

\begin{Lem} \label{L:embedding}
For each $ p = \left( H, \pi, \mathcal{N}, T \right) \in \mathbb{P}$,
there exists  
$p^{+} = \left( H, \pi^{+}, \mathcal{N}^{+}, T \right) \leq p$
such that $\pi^{+} : H \to \Aut \mathcal{N}^{+}$ is an embedding.
\end{Lem}

\begin{proof}
Let $\mathcal{N}^{+}$ be the structure for the language
$L( \mathcal{N})$ such that
\begin{enumerate}
\item[(a)] the universe of $\mathcal{N}^{+}$ is the disjoint union
$\mathcal{N} \sqcup H$;
\item[(b)] for each relation $R \in L( \mathcal{N})$,
$R^{\mathcal{N}^{+}} = R^{\mathcal{N}}$.
\end{enumerate}
In particular, if $x \in \mathcal{N}^{+} \smallsetminus \mathcal{N}$,
then $x$ only realises the trivial restricted atomic type $\emptyset$ over
$\mathcal{N}$. Hence none of the restricted atomic types in $T$ is realised in
$\mathcal{N}^{+}$. Let
$\pi^{+} : H \to \Aut \mathcal{N}^{+}$ be the embedding 
such that for each $h \in H$,
\begin{enumerate}
\item[(i)] $\pi^{+}(h)(x) = \pi(h)(x)$ for all $x \in \mathcal{N}$; and
\item[(ii)] $\pi^{+}(h)(x) = hx$ for all $x \in H$.
\end{enumerate}
Then $p^{+} = \left( H, \pi^{+}, \mathcal{N}^{+}, T \right) \leq p$.
\end{proof}

There is a slight inaccuracy in the proof of Lemma \ref{L:embedding},
as the universe of $\mathcal{N}^{+}$ should really be an ordinal
$\delta < \kappa$. However, the proof can easily be repaired: simply
replace $\mathcal{N}^{+}$ by a suitable isomorphic structure.
Similar remarks apply to the proofs of Lemmas \ref{L:domain} and
\ref{L:onto}.

\begin{Lem} \label{L:domain}
For each $a \in G$, 
\[
D_{a} = \{ \left( H, \pi, \mathcal{N}, T \right) \mid a \in H \}
\]
is a dense subset of $\mathbb{P}$.
\end{Lem}

\begin{proof}
Let $a \in G$ and $p = \left( H, \pi, \mathcal{N}, T \right) \in
\mathbb{P}$. We can suppose that $a \notin H$. Let
$H^{+} = \langle H, a \rangle$. Let
$C = \{ g_{i} \mid i \in I \}$ be a set of left coset representatives
for $H$ in $H^{+}$, chosen so that $1 \in C$. Let 
$\mathcal{N}^{+}$ be the structure for the language
$L( \mathcal{N})$ such that
\begin{enumerate}
\item[(a)] the universe of $\mathcal{N}^{+}$ is the
cartesian product $C \times \mathcal{N}$; and
\item[(b)] for each relation $R \in 
L(\mathcal{N})$,
\[
\left( (g_{i},x) , (g_{j},y) \right) \in R^{\mathcal{N}^{+}}
\text{ if{f} } i=j \text{ and } \left( x,y \right) \in R^{\mathcal{N}}.
\]
\end{enumerate}
By identifying each $x \in \mathcal{N}$ with the element
$( 1, x) \in \mathcal{N}^{+}$, we can regard 
$\mathcal{N}$ as a substructure of $\mathcal{N}^{+}$.
Once again, each element 
$(g_{i},x) \in \mathcal{N}^{+} \smallsetminus \mathcal{N}$ only
realises the trivial restricted atomic type $\emptyset$ over
$\mathcal{N}$; and hence none of the restricted atomic types in $T$ is realised
in $\mathcal{N}^{+}$.

Define an action of $H^{+}$ on $\mathcal{N}^{+}$ as follows.
If $g \in H^{+}$ and $(g_{i},x) \in \mathcal{N}^{+}$, then
\[
g (g_{i},x) = ( g_{j}, \pi(h)(x)),
\]
where $j \in I$ and $h \in H$ are such that
$g g_{i} = g_{j}h$. It is easily checked that this action yields a
homomorphism $\pi^{+} : H^{+} \to \Aut \mathcal{N}^{+}$;
and that $\left( H^{+}, \pi^{+}, \mathcal{N}^{+}, T \right)
\leq p$.
\end{proof}

\begin{Lem} \label{L:size}
For each $\alpha < \kappa$,
\[
E_{\alpha} = \{\left( H, \pi, \mathcal{N}, T \right) 
\mid \alpha \in \mathcal{N} \}
\]
is a dense subset of $\mathbb{P}$.
\end{Lem}

\begin{proof}
Left to the reader.
\end{proof}

Let $F \subseteq \mathbb{P}$ be a generic filter, and let
$V^{\mathbb{P}} = V[F]$ be the corresponding generic
extension. Working within $V^{\mathbb{P}}$, define
\[
\mathcal{M} = \bigcup \{ \mathcal{N} \mid
\text{ There exists } p = \left( H, \pi, \mathcal{N}, T \right) \in F \}
\]
and \[
\pi^{*} = \bigcup \{ \pi \mid
\text{ There exists } p = \left( H, \pi, \mathcal{N}, T \right) \in F \}.
\]
Then $\mathcal{M}$ is a structure for $L$ of cardinality $\kappa$,
and $\pi^{*}$ is an embedding of $G$ into 
$\Aut \mathcal{M}$. So the following lemma completes the proof
of Theorem \ref{T:main2}.

\begin{Lem} \label{L:onto}
$\pi^{*} : G \to \Aut \mathcal{M}$ is a surjective homomorphism.
\end{Lem}

\begin{proof}
Suppose that $g \in \Aut \mathcal{M} \smallsetminus \pi^{*}[G]$. Let
$\widetilde{\mathcal{M}}$, $\widetilde{\pi}$ be the canonical
$\mathbb{P}$-names for $\mathcal{M}$ and $\pi^{*}$; and let
$\widetilde{g}$ be a $\mathbb{P}$-name for $g$. Then there exists a
condition $p \in F$ such that
\[
p \Vdash \widetilde{g} \in \Aut \widetilde{\mathcal{M}} \text{ and }
\widetilde{g} \neq \widetilde{\pi}(h) \text{ for all } h \in G.
\]
Using the fact that $\mathbb{P}$ is $\kappa$-closed, we can inductively
construct a descending sequence of conditions
$ p_{m} = \left( H_{m}, \pi_{m}, \mathcal{N}_{m}, T_{m} \right)$
for $m \in \omega$ such that the following hold.
\begin{enumerate}
\item[(a)] $p_{0} = p$.
\item[(b)] For all $x \in \mathcal{N}_{m}$, there exists 
$y \in \mathcal{N}_{m+1}$
such that $p_{m+1} \Vdash \widetilde{g}(x) = y$.
\item[(c)] For all $h \in H_{m}$, there exists 
$z \in \mathcal{N}_{m+1}$ such that
$p_{m+1} \Vdash \widetilde{g}(z) \neq \pi_{m+1}(h)(z)$.
\end{enumerate}
Let $q = \left( H, \pi, \mathcal{N}, T \right)$ be the greatest lower
bound of $\{ p_{m} \mid m \in \omega \}$ in $\mathbb{P}$.
Then $q \leq p$, and there exists
$g^{*} \in \Aut \mathcal{N} \smallsetminus \pi [H]$ such that
$q \Vdash \widetilde{g} \res \mathcal{N} = g^{*}$. Let
$\mathcal{N}^{+}$ be the structure defined as follows.
\begin{enumerate}
\item[(1)] The universe of $\mathcal{N}^{+}$ is the disjoint
union $\mathcal{N} \sqcup H$.
\item[(2)] For each relation $R \in L(\mathcal{N})$,
$R^{\mathcal{N}^{+}} = R^{\mathcal{N}}$.
\item[(3)] For each $x \in \mathcal{N}$, let $R_{x} \in 
L \smallsetminus L(\mathcal{N})$ 
be a new binary relation symbol. Then we set
$(h,y) \in R_{x}^{\mathcal{N}^{+}}$ if{f} $h \in H$,
$y \in \mathcal{N}$ and $\pi(h)(x) = y$.
\end{enumerate}
Once again, it is clear that none of the restricted atomic types in $T$ is realised in
$\mathcal{N}^{+}$. Let
$\pi^{+} : H \to \Sym \mathcal{N}^{+}$ be the embedding
such that
\begin{enumerate}
\item[(i)] $\pi^{+}(h)(x) = \pi(h)(x)$ for all $x \in \mathcal{N}$; and
\item[(ii)] $\pi^{+}(h)(x) = hx$ for all $x \in H$.
\end{enumerate}
Then it is easily checked that
$\pi^{+}[H] \leqslant \Aut \mathcal{N}^{+}$. Finally let $t$ be the 
restricted atomic type defined by
\[
t = \{ R_{x}(v, g^{*}(x) \mid x \in \mathcal{N} \};
\]
and let $T^{+} = T \cup \{ t \}$.

\begin{Claim} \label{C:onto}
$t$ is omitted in $\mathcal{N}^{+}$.
\end{Claim}

\begin{proof}[Proof of Claim \ref{C:onto}]
If $z \in \mathcal{N}^{+}$ realises $t$, then $z = h \in H$.
And if $x \in \mathcal{N}$, then 
$\mathcal{N}^{+} \vDash R_{x}(h, g^{*}(x))$, and so
$\pi(h)(x) = g^{*}(x)$. But this contradicts the fact that
$g^{*} \in \Aut \mathcal{N} \smallsetminus \pi [H]$.
\end{proof}

Thus $q^{+} = \left( H, \pi^{+}, \mathcal{N}^{+}, T^{+} \right)
\in \mathbb{P}$. To simplify notation, suppose that
$q^{+} \in F$; so that $g^{*} \subseteq g$. Then for each
$x \in \mathcal{N}$, we have that
$\mathcal{M} \vDash R_{x}(1,x)$, and hence
$\mathcal{M} \vDash R_{x}(g(1), g^{*}(x))$. But this means
that $g(1) \in \mathcal{M}$ realises $t$, which is the final contradiction.
\end{proof}

\section{$\tau_{\kappa}$ is increasing} \label{S:increasing}

In this section, we will complete the proof of Theorem \ref{T:best}. So
let $V \vDash GCH$ and let $\kappa$, $\lambda  \in V$ be uncountable cardinals
such that $\kappa < \cf (\lambda)$. Let $\alpha$ be any ordinal such
that $\alpha < \lambda^{+}$. It is well-known that there exists a
centreless group $T$ of cardinality $\kappa$ such that
$\Aut T = \Inn T$.
(For example, we can take $T = PGL(2,K)$, where
$K$ is a rigid field of cardinality $\kappa$. The existence of such
a field follows from Lemma \ref{L:field}.) Thus we can assume that
$\alpha \geq 1$. First consider the case when $\kappa$ is
a regular cardinal. Let $\mathbb{Q}$ be the notion of forcing which is
given by Theorem \ref{T:just}. Then the following statements are true
in the corresponding generic extension $M = V^{\mathbb{Q}}$.
\begin{enumerate}
\item[(a)] $\kappa^{< \kappa} = \kappa$.
\item[(b)] $2^{\kappa} = \lambda$.
\item[(c)] There exist groups $H < G < \Sym (\kappa)$ such that
$|H| = \kappa$ and the normaliser tower of $H$ in $G$ terminates in
exactly $\alpha$ steps.
\end{enumerate}
Let $\mathbb{P} \in M$ be the notion of forcing, given by Theorem
\ref{T:main2}, which adjoins a graph $\Gamma$ of cardinality
$\kappa$ such that $G \simeq \Aut \Gamma$. Then, applying Lemma
\ref{L:reduction}, we find that the following statements
are true in $M^{\mathbb{P}}$.
\begin{enumerate}
\item[(1)] $2^{\kappa} = \lambda$.
\item[(2)] There exists a centreless group $T$ of cardinality
$\kappa$ such that $\tau (T) = \alpha$. 
\end{enumerate}

Next suppose that $\kappa$ is a singular cardinal. By the above argument,
there is a generic extension $V^{\mathbb{P} \ast \mathbb{Q}}$ in which
the following statements are true.
\begin{enumerate}
\item[(i)] $2^{\omega_{1}} = 2^{\kappa} = \lambda$.
\item[(ii)] There exists a centreless group $T$ of cardinality
$\omega_{1}$ such that $\tau (T) = \alpha$. 
\end{enumerate}
Let $G = T \times \Alt (\kappa)$. Clearly $|G| = \kappa$;
and the following theorem implies that 
$\tau (G) = \tau (T) = \alpha$. 
This completes the proof of Theorem \ref{T:best}.

\begin{Thm} \label{T:increasing}
Suppose that $\omega \leq \theta < \kappa$. If $H$ is a centreless group of
cardinality $\theta$ such that $\tau (H) \geq 1$, then
$\tau (H \times \Alt (\kappa)) = \tau (H)$.
\end{Thm}

\begin{Cor} \label{C:increasing}
If $\omega \leq \theta < \kappa$, then
$\tau_{\theta} \leq \tau_{\kappa}$.
\end{Cor}

The remainder of this section will be devoted to the proof
of Theorem \ref{T:increasing}. Let $G = H \times \Alt (\kappa)$; and
let $H_{\beta}$, $G_{\beta}$ be the $\beta^{th}$ groups in the automorphism
towers of $H$, $G$ respectively. We will eventually prove by induction
that $G_{\beta} = H_{\beta} \times \Sym (\kappa)$ for all $\beta \geq 1$.
To accomplish this, we need to keep track of $\varphi [\Alt (\kappa) ]$
for each automorphism $\varphi$ of $G_{\beta}$. The next lemma shows
that for all $\varphi \in \Aut G_{\beta}$, either 
$\varphi [\Alt (\kappa) ] \leqslant H_{\beta}$ or 
$\varphi [\Alt (\kappa) ] \leqslant \Sym (\kappa)$. The main point will
be to eliminate the possibility that
$\varphi [\Alt (\kappa) ] \leqslant H_{\beta}$. This is straightforward
when $\beta$ is a successor ordinal. To deal with the case when $\beta$
is a limit ordinal, we will make use of the result that
$\Alt (\kappa)$ is strictly simple.

\begin{Lem} \label{L:increasing1}
Suppose that $A$ is a simple nonabelian normal subgroup of the direct product
$H \times S$. Then either $A \leqslant H$ or
$A \leqslant S$.
\end{Lem}

\begin{proof}
Let $1 \neq g = xy \in A$, where $x \in H$ and $y \in S$. If $y=1$,
then the conjugacy class $g^{A} = x^{A}$ is contained in $H$, and
so $A = \langle g^{A} \rangle \leqslant H$. So suppose that
$y \neq 1$. Let $\pi : H \times S \to S$ be the canonical projection
map. Then $1 \neq y \in \pi [A] \leqslant S$ and 
$\pi [A] \simeq A$. Hence there exists an element
$z \in \pi [A] \leqslant S$ such that $z y z^{-1} \neq y$. Since
$A \trianglelefteq H \times S$, it follows that
\[
1 \neq zyz^{-1}y^{-1} = zxyz^{-1}y^{-1}x^{-1} = zgz^{-1}g^{-1}
\in A \cap S.
\]
Arguing as above, we now obtain that
$A \leqslant S$.
\end{proof}

\begin{Def} \label{D:ascendant}
Let $H$ be a subgroup of the group $G$. Then $H$ is said to be an
{\em ascendant subgroup\/} of $G$ if there exist an ordinal $\gamma$ and a set
of subgroups $\{ H_{\beta} \mid \beta \leq \gamma \}$ such that the following
conditions are satisfied.
\begin{enumerate}
\item[(a)] $H_{0} = H$ and $H_{\gamma} = G$.
\item[(b)] If $\beta < \gamma$, then 
$H_{\beta} \trianglelefteq H_{\beta +1}$.
\item[(c)] If $\delta$ is a limit ordinal such that $\delta \leq \gamma$, then
$H_{\delta} = \bigcup_{\beta < \delta}H_{\beta}$. 
\end{enumerate}
\end{Def}

\begin{Def} \label{D:strictly}
A group $A$ is {\em strictly simple\/} if it has no nontrivial
proper ascendant subgroups.
\end{Def}

\begin{Thm}[Macpherson and Neumann \cite{mn}] \label{T:strictly}
For each $\kappa \geq \omega$, the alternating group $\Alt (\kappa)$
is strictly simple.
\end{Thm}
\begin{flushright}
$\square$
\end{flushright}

\begin{Lem} \label{L:increasing2}
Suppose that $H$ is a centreless group. Let $\tau = \tau (H)$ and let
$H_{\tau}$ be the $\tau^{th}$ group in the 
automorphism tower of $H$. If $A$ is a strictly simple normal subgroup
of $H_{\tau}$, then $A \leqslant H_{0}$.
\end{Lem}

\begin{proof}
Let $\beta \leq \tau$ be the least ordinal such that
$A \leqslant H_{\beta}$. First suppose that $\beta$ is a
limit ordinal. Then
$A = \bigcup_{\gamma < \beta }\left(A \cap H_{\gamma}\right)$,
and $A \cap H_{\gamma} \trianglelefteq A \cap H_{\gamma +1}$ for all
$\gamma < \beta$. Consequently, if $\gamma < \beta$ is the least
ordinal such that $A \cap H_{\gamma} \neq 1$, then $A \cap H_{\gamma}$
is a nontrivial ascendant subgroup of $A$. But this contradicts the
assumption that $A$ is strictly simple.

Next suppose that $\beta = \gamma +1$ is a successor ordinal. Since
$A \cap H_{\gamma}$ is a proper normal subgroup of $A$, it follows
that $A \cap H_{\gamma} =1$. Now notice that
$A \leqslant H_{\gamma +1} \leqslant N_{H_{\tau}}(H_{\gamma})$ and 
$H_{\gamma} \leqslant H_{\tau} = N_{H_{\tau}}(A)$, 
This implies that
$\left[ A, H_{\gamma} \right] \leqslant A \cap H_{\gamma} =1$.
(For example, see Lemma 1.1.3 of Suzuki \cite{suz}.)
But then
$A \leqslant C_{H_{\gamma +1}}(H_{\gamma}) =1$, which is a
contradiction. The only remaining possibility is that
$\beta =0$.
\end{proof}

We will also make use of the well-known results that
$\Aut ( \Alt (\kappa)) = \Sym (\kappa)$ and
$\Aut ( \Sym (\kappa)) = \Sym (\kappa)$. 
(For example, see Theorem 11.4.8 of Scott \cite{scott}.)

\begin{proof}[Proof of Theorem \ref{T:increasing}]
Let $\tau = \tau (H)$; and let
\[
H = H_{0} \trianglelefteq H_{1} \trianglelefteq  \dots
H_{\beta} \trianglelefteq H_{\beta +1} \trianglelefteq \dots
H_{\tau} = H_{\tau +1} = \dots
\]
be the automorphism tower of of $H$. Let $G = H \times \Alt (\kappa)$;
and let
\[
G = G_{0} \trianglelefteq G_{1} \trianglelefteq G_{2} \trianglelefteq  \dots
G_{\beta} \trianglelefteq G_{\beta +1} \trianglelefteq \dots
\]
be the automorphism tower of $G$. We will prove by induction 
on $\beta \geq 1$ that
$G_{\beta} = H_{\beta} \times \Sym (\kappa)$.

First consider the case when $\beta =1$. Let $\varphi \in \Aut G$ be
any automorphism. Then $\varphi [\Alt (\kappa)]$ is a simple
nonabelian normal subgroup of the direct product
$G = H \times \Alt (\kappa)$. By Lemma \ref{L:increasing1}, either
$\varphi [\Alt (\kappa)] \leqslant H$ or
$\varphi [\Alt (\kappa)] \leqslant \Alt (\kappa)$. Since
$\left| \varphi [\Alt (\kappa)] \right| = \kappa
> \theta = |H|$, it follows that
$\varphi [\Alt (\kappa)] \leqslant \Alt (\kappa)$. As
$\varphi [\Alt (\kappa)]$ is a normal subgroup of $G$, we must
have that $\varphi [\Alt (\kappa)] = \Alt (\kappa)$. Note that
$C_{G}( \Alt (\kappa)) = H$. Hence we must also have that
$\varphi [H] =H$. It follows that
\[
G_{1} = \Aut G = \Aut H \times \Aut (\Alt (\kappa)) = H_{1} \times
\Sym (\kappa).
\]

Next suppose that $\beta = \gamma +1$ and that
$G_{\gamma} = H_{\gamma} \times \Sym (\kappa)$. Let 
$\varphi \in \Aut G_{\gamma}$ be any automorphism. By Lemma
\ref{L:increasing1},  either
$\varphi [\Alt (\kappa)] \leqslant H_{\gamma}$ or
$\varphi [\Alt (\kappa)] \leqslant \Sym (\kappa)$. As $\Alt (\kappa)$
is a strictly simple group, Lemma \ref{L:increasing2} implies that
$\varphi [\Alt (\kappa)] \leqslant \Sym (\kappa)$. Since
$\varphi [\Alt (\kappa)]$ is a simple normal subgroup of $G_{\gamma}$,
it follows that $\varphi [\Alt (\kappa)] = \Alt (\kappa)$. Using the
facts that 
$C_{G_{\gamma}}( \Alt (\kappa)) = H_{\gamma}$ and
$C_{G_{\gamma}}(H_{\gamma}) = \Sym(\kappa)$, we now see that
$\varphi [H_{\gamma}] = H_{\gamma}$ and 
$\varphi [\Sym (\kappa)] = \Sym (\kappa)$. Hence
\[
G_{\gamma +1} = \Aut G_{\gamma} 
= \Aut H_{\gamma} \times \Aut (\Sym (\kappa)) 
= H_{\gamma +1} \times \Sym (\kappa).
\]

Finally no difficulties arise when $\beta$ is a limit ordinal.
\end{proof}


\begin{thebibliography}{99}

\bibitem{b}
N. G. De Bruijn, {\em Embedding theorems for infinite groups\/},
Indag. Math. {\bf 19} (1957), 560--569.

\bibitem{fk}
E. Fried and J. Koll\'{a}r, {\em Automorphism groups of fields\/},
in Universal Algebra (E. T. Schmidt et al., eds.), Coloq. Math.
Soc. Janos Boyali, vol {\bf 24}, 1981, pp. 293--304.

\bibitem{h}
W. Hodges, {\em Model theory\/}, Encyclopedia of Mathematics
and its Applications {\bf 42}, Cambridge University Press, 1993.

\bibitem{hu}
J. A. Hulse, {\em Automorphism towers of polycyclic groups\/},
J. Algebra {\bf 16} (1970), 347--398.

\bibitem{ku}
D. W. Kueker, {\em Definability, automorphisms and infinitary
languages\/}, in The Syntax and Semantics of Infinitary Languages
(ed. J. Barwise), Lecture Notes in Math. {\bf 72}, Springer, Berlin,
1968, pp. 152--165.

\bibitem{k}
K. Kunen, {\em Set Theory\/}, North-Holland, Amsterdam, 1980.

\bibitem{mn}
H. D. Macpherson and P. M. Neumann, {\em Subgroups of infinite
symmetric groups\/}, J. London Math. Soc. (2) {\bf 42} (1990),
64--84.

\bibitem{rr}
A. Rae and J. E. Roseblade, {\em Automorphism towers of extremal
groups\/}, Math. Z. {\bf 117} (1970), 70--75.

\bibitem{scott}
W. R. Scott, {\em Group Theory\/}, Prentice-Hall,
Englewood Cliffs, New Jersey, 1964.

\bibitem{s}
S. Shelah, {\em First order theory of permutation groups\/}, Israel J. Math.
{\bf 14} (1973), 149--162.

\bibitem{so}
S. Solecki, {\em Polish group topologies\/}, preprint (1997).

\bibitem{suz}
M. Suzuki, {\em Group Theory II\/}, Grundlehren der mathematischen
Wissenschaften {\bf 248}, Springer-Verlag 1986.

\bibitem{t1}
S. Thomas, {\em The automorphism tower problem\/}, Proc. Amer. Math.
Soc. {\bf 95} (1985), 166--168.

\bibitem{t2}
S. Thomas, {\em The automorphism tower problem II\/}, to appear in
Israel J. Math.

\bibitem{w}
H. Wielandt, {\em Eine Verallgemeinerung der invarianten Untergruppen\/},
Math. Z. {\bf 45} (1939), 209--244.


\end{thebibliography}
\end{document}